\documentclass[12pt]{article}

\usepackage{cite,amsfonts,amssymb}
\newcommand{\blackbox}{\vrule Depth0pt height7pt width7pt}
\newcommand{\bgproof}{\noindent {\bf Proof} \hspace{2mm}}
\newcommand{\edproof}{\hfill \blackbox \vspace{3mm}}

\newtheorem{prop}{Proposition}

\newtheorem{lem}{Lemma}

\newtheorem{theor}{Theorem}

\newtheorem{cor}{Corollary}

\newtheorem{rem}{Remark}

\newtheorem{defin}{Definition}

\date{}

\title{\bf Stochastic Volterra equations driven by cylindrical Wiener process}

\author{{\large\sf Anna Karczewska and Carlos Lizama}\\[2mm]
  \normalsize\it
 Department of Mathematics,
 University of Zielona G\'ora\\ \normalsize\it
 ul. Szafrana 4a, 65-246 Zielona G\'ora, Poland,~
 e-mail: A.Karczewska@im.uz.zgora.pl\\[2mm] \normalsize\it
 Universidad de Santiago de Chile, Departamento de
Matem\'atica, Facultad de Ciencias\\ \normalsize\it Casilla
307-Correo 2, Santiago, Chile,~ e-mail: clizama@lauca.usach.cl \\
}

\begin{document}

\maketitle

\def\thefootnote{}

\vspace{-10mm}
\begin{abstract}

In this paper, stochastic Volterra equations driven by cylindrical
Wiener process in Hilbert space are investigated.
Sufficient conditions for existence of strong solutions are given.
The key role is played by convergence of $\alpha$-times resolvent
families.

\end{abstract}

\section{Introduction}\label{sSW1}
Let $H$ be a separable Hilbert space with a norm $|\cdot|_H$ and
$A$ be a closed linear unbounded operator with dense domain
$D(A)\subset H$ equipped with the graph norm $|\cdot|_{D(A)}$. The
purpose of this paper is to study the existence of strong
solutions for a class of  stochastic Volterra equations of the
form
\begin{equation}{\label{eq4}}
X(t) = X_0 + \int_0^t a(t-\tau) AX(\tau)d\tau + \int_0^t
\Psi(\tau)\, dW(\tau), \quad t\geq 0,
\end{equation}
where $ a(t) =\displaystyle \frac{t^{\alpha-1}}{\Gamma(\alpha)},~
\alpha>0,$ and $W, \Psi$ are appropriate stochastic processes. It
is well known that there are several situations that can be
modeled by stochastic Volterra equations (see e.g. \cite[Section
3.4 ]{Ho-Ok-Ub-Zh} and references therein). We note that
stochastic Volterra equations driven by white noise have been
studied in \cite{Cl-Pr96} among other authors. A similar equation
and very related to our case appears first studied in
\cite{Bo-Tu}. Here we are interested in the study of strong
solutions when equation (\ref{eq4}) is driven by a
cylindrical Wiener process $W$.

 When $a(t)$ is a completely positive function,
sufficient conditions for existence of strong solutions for
(\ref{eq4}) were obtained in \cite{KL2}. This was done using a
method which involves the use of a resolvent family associated to
the  deterministic version of equation (\ref{eq4}):
\begin{equation} \label{eq3a}
u(t) = \int_0^t a(t-\tau)\,Au(\tau)d\tau + f(t), \quad t\geq 0,
\end{equation}
where $f$ is an $H$-valued function.

However, there are two kinds of problems that arise when we study
(\ref{eq4}). On the one hand, the kernels
$\frac{t^{\alpha-1}}{\Gamma(\alpha)} $ are $\alpha$-regular and
$\frac{\alpha \pi}{2}$ -sectorial but not completely positive
functions for $\alpha >1$, so e.g. the results in \cite{KL2}
cannot be used directly for $\alpha >1$. On the other hand, for
$\alpha\in (0,1)$, we have a singularity of the kernel in $t=0.$
This fact strongly suggests the use of
 $\alpha$-times resolvent families associated to
equation (\ref{eq3a}). This new tool appears carefully studied in
\cite{Ba} as well as their relationship with fractional
derivatives. For convenience of the reader, we
  provide below the main results on $\alpha$-times resolvent
families to be used in this paper.

Our second main ingredient to obtain strong solutions of
(\ref{eq4}) relies on approximation of $\alpha$-times resolvent
families. This kind of result was very recently formulated by Li
and Zheng \cite{LiZh}. It enables us to prove a key result on
convergence of $\alpha$-times resolvent families (see Theorem \ref{th2}
below). Then we can follow the methods employed in \cite{KL2} to obtain
existence of strong solution for the stochastic equation
(\ref{eq4}) (see Theorem \ref{coSW4}).

Our plan for the paper is the following. In section 2 we formulate
the deterministic results which will play the key role for the
paper. Section 3 is devoted to weak and mild solutions
while in section 4 we provide strong
solution to (\ref{eq4}). More precisely, we give sufficient condition
for a stochastic convolution to be a strong solution to (\ref{eq4}).

\section{Convergence of $\alpha$-times resolvent families}\label{Sconv}

In this section we formulate the main deterministic results on
convergence of resolvents. We denote
$$
g_{\alpha}(t) = \frac{t^{\alpha -1}}{\Gamma (\alpha)}, \quad
\alpha > 0, \quad t >0,
$$
where $\Gamma$ is the gamma function.

By $S_{\alpha}(t),~t\geq 0$, we denote the family of
$\alpha$-times resolvent families corresponding to the Volterra
equation (\ref{eq3a}), if it exists, and defined as follows.

\begin{defin}\label{def1} (see \cite{Ba})\\
A family $(S_{\alpha}(t))_{t\geq 0}$ of bounded linear operators
in a Banach space $B$ is called {\tt $\alpha$-times resolvent
family} for (\ref{eq3a}) if the following conditions are
satisfied:
\begin{enumerate}
\item $S_{\alpha}(t)$ is strongly continuous on $\mathbb{R}_+$ and
$S_{\alpha}(0)=I$;
\item $S_{\alpha}(t)$ commutes with the operator $A$, that is,
$S_{\alpha}(t)(D(A))\subset D(A)$ and
$AS_{\alpha}(t)x=S_{\alpha}(t)Ax$ for all $x\in D(A)$ and $t\geq
0$;
\item the following {\tt resolvent equation} holds
\begin{equation} \label{eq4a}
S_{\alpha}(t)x = x + \int_0^t g_{\alpha}(t-\tau)
AS_{\alpha}(\tau)x d\tau
\end{equation}
for all $x\in D(A),~t\geq 0$.
\end{enumerate}
\end{defin}

Necessary and sufficient conditions for existence of the
$\alpha$-times resolvent family have been studied in \cite{Ba}.
Observe that the $\alpha$-times resolvent family corresponds to a
$C_0$-semigroup in case $\alpha =1$ and a cosine family in case
$\alpha=2.$ In consequence, when $1 < \alpha < 2$ such resolvent
families interpolate $C_0$-semigroups and cosine functions. In
particular, for $A = \Delta$, the integrodifferential equation
corresponding to such resolvent family interpolates the heat
equation and the wave equation (see \cite{Fu}).

\begin{defin}{\label{def2}}
An $\alpha$-times resolvent family $(S_{\alpha}(t))_{t\geq 0}$ is
called \mbox{\tt exponentially bounded} if there are constants $M
\geq 1$ and $\omega \geq 0$ such that
\begin{equation}{\label{eq5}}
\| S_{\alpha}(t) \| \leq Me^{\omega t}, \quad t\geq 0.
\end{equation}
\end{defin}

If there is the $\alpha$-times resolvent family
$(S_{\alpha}(t))_{t\geq 0}$ for $A$ and satisfying (\ref{eq5}), we
write $A \in \mathcal{C}^{\alpha}(M,\omega). $ Also, set $
\mathcal{C}^{\alpha}(\omega) :=\cup_{M\geq 1}
\mathcal{C}^{\alpha}(M,\omega)$ and $ \mathcal{C}^{\alpha}
:=\cup_{\omega \geq 0} \mathcal{C}^{\alpha}(\omega).$

\begin{rem}\label{rem1}{\rm
It was proved by Bazhlekova \cite[Theorem 2.6]{Ba} that if $ A\in
\mathcal{C}^{\alpha}$ for some $\alpha > 2,$ then $A$ is bounded.}
\end{rem}

The following subordination principle is very important in the
theory of $\alpha$-times resolvent families (see \cite[Theorem
3.1]{Ba}).

\begin{theor}\label{th1}
Let $0 < \alpha < \beta \leq 2, \gamma = \alpha /\beta, \omega\geq
0.$ If $A \in \mathcal{C}^{\beta}(\omega)$ then $A \in
\mathcal{C}^{\alpha}(\omega^{1/\gamma})$ and the following
representation holds
\begin{equation}\label{eq6}
S_{\alpha}(t)x = \int_0^{\infty} \varphi_{t,\gamma}(s)
S_{\beta}(s)xds, \quad t>0,
\end{equation}
where $\varphi_{t,\gamma}(s) := t^{-\gamma}
\Phi_{\gamma}(st^{-\gamma})$ and $ \Phi_{\gamma}(z)$ is the Wright
function defined as
\begin{equation}\label{eq7}
\Phi_{\gamma}(z):= \sum_{n=0}^{\infty} \frac{(-z)^n}
{n!\,\Gamma(-\gamma n + 1 - \gamma)}, \quad 0 < \gamma < 1.
\end{equation}
\end{theor}

\begin{rem}\label{rem2}{\rm
(i) We recall that the Laplace transform of the Wright function
corresponds to $E_{\gamma}(-z)$ where $E_{\gamma}$ denotes the
Mittag-Leffler function. In particular, $\Phi_{\gamma}(z)$ is a
probability density function.

(ii) Also we recall from \cite[(2.9)]{Ba} that the continuity in
$t\geq 0$ of the Mittag-Leffler function together with the
asymptotic behavior of it, imply that for $\omega \geq 0$ there
exists a constant $C>0$ such that
\begin{equation}\label{eq8}
E_{\alpha}(\omega t^{\alpha}) \leq C e^{{\omega^{1/\alpha}} t},
\quad t \geq 0,\,\, \alpha \in (0,2).
\end{equation}
}\end{rem}

As we have already written, in this paper the results concerning
convergence of $\alpha$-times resolvent families in a Banach space
$B$ will play the key role. Using a very recent result due to Li
and Zheng \cite{LiZh} we are able to prove the following theorem.

\begin{theor} \label{th2}
Let $A$ be the generator of a $C_0$-semigroup $(T(t))_{t\geq 0}$
in a Banach space $B$ such that
\begin{equation}\label{eq9}
\|T(t) \| \leq Me^{\omega t}, \quad t \geq 0.
\end{equation}
Then, for each $ 0 < \alpha < 1$ we have $A \in
\mathcal{C}^{\alpha}(M,\omega^{1/\alpha}). $ Moreover, there exist
bounded operators $A_n$ and $\alpha$-times resolvent families
$S_{\alpha,n}(t)$ for $A_n$ satisfying $ ||S_{\alpha,n}(t) || \leq
MCe^{(2\omega)^{1/\alpha} t},$ for all $t\geq 0,~n\in \mathbb{N},$
and
\begin{equation} \label{eq10}
S_{\alpha,n}(t)x \to S_{\alpha}(t)x \quad \mbox{as} \quad n\to
+\infty
\end{equation}
for all $x \in B,\; t\geq 0.$ Moreover, the convergence is uniform
in $t$ on every compact subset of $ \mathbb{R}_+$.
\end{theor}

\bgproof Since $A$ is the generator of a $C_0$ semigroup
satisfying (\ref{eq9}), we have $A\in \mathcal{C}^1(\omega)$.
Hence, the first assertion follows directly from Theorem
\ref{th1}, that is, for each $ 0<\alpha < 1$ there is an
$\alpha$-times resolvent family $(S_{\alpha}(t))_{t\geq 0}$ for
$A$ given by
\begin{equation}\label{eq11}
S_{\alpha}(t)x = \int_0^{\infty} \varphi_{t,\alpha}(s) T(s)xds,
\quad t>0.
\end{equation}

Since $A$ generates a $C_0$-semigroup, the resolvent set $\rho(A)
$ of $A$ contains the ray $ [w,\infty)$ and
$$
||R(\lambda,A)^k || \leq \frac{M}{(\lambda -w)^k } \qquad
\mbox{for } \lambda > w, \qquad k\in \mathbb{N}.
$$

Define
\begin{equation} \label{eq12}
A_n := n AR(n,A) = n^2 R(n,A) - nI, \qquad n> w,
\end{equation}
the {\it Yosida approximation} of $A$.

Then
\begin{eqnarray*}
||e^{t A_n} || &=& e^{-nt} || e^{n^2 R(n,A)t} || \leq
e^{-nt} \sum_{k=0}^{\infty} \frac{n^{2k} t^k}{k!} ||R(n,A)^k|| \\
&\leq& M e^{(-n + \frac{n^2}{n-w})t} = M e^{ \frac{nwt}{ n-w}}.
\end{eqnarray*}
Hence, for $n > 2w$ we obtain
\begin{equation}{\label{eq13}}
|| e^{A_n t} || \leq M e^{2wt}.
\end{equation}
Next, since each $A_n$ is bounded, it follows also from Theorem
\ref{th1} that for each $0< \alpha < 1$ there exists an
$\alpha$-times resolvent family $ (S_{\alpha,n}(t))_{t\geq 0}$ for
$A_n$ given as

\begin{equation}{\label{eq14}}
S_{\alpha,n}(t) = \int_0^{\infty} \varphi_{t,\alpha}(s)
e^{sA_n}ds, \quad t>0.
\end{equation}

By (\ref{eq13}) and Remark \ref{rem2}(i) it follows that

\begin{eqnarray*}
\| S_{\alpha,n}(t) \| &\leq & \int_0^{\infty}
\varphi_{t,\alpha}(s) \| e^{s A_n} \| ds \\ &\leq& M
\int_0^{\infty} \varphi_{t,\alpha}(s) e^{2\omega s} ds = M
\int_0^{\infty} \Phi_{\alpha}(\tau) e^{2\omega t^{\alpha} \tau }
d\tau= M E_{\alpha}(2 \omega t^{\alpha}), \quad t \geq 0.
\end{eqnarray*}

This together with Remark \ref{rem2}(ii), gives

\begin{equation}
\| S_{\alpha,n}(t) \| \leq MCe^{(2\omega)^{1/\alpha}t}, \quad t
\geq 0.
\end{equation}

Now, we recall the fact that $ R(\lambda,A_n)x \to R(\lambda,A)x $
as $ n\to \infty$ for all $\lambda $ sufficiently large (see e.g.
\cite[Lemma~7.3]{Pa}), so we can conclude from \cite[Theorem
4.2]{LiZh} that
\begin{equation}
S_{\alpha,n}(t)x \to S_{\alpha}(t)x \quad \mbox{as} \quad n\to
+\infty
\end{equation}
for all $x \in B,$ uniformly for $t$ on every compact subset of
$\mathbb{R}_+\,$. \edproof

An analogous result can be proved in the case when $A$ is the
generator of a strongly continuous cosine family.

\begin{theor} \label{th3a}
Let $A$ be the generator of a $C_0$-cosine family $(T(t))_{t\geq
0}$ in a Banach space $B$. Then, for each $0<\alpha<2$ we have $A
\in \mathcal{C}^{\alpha}(M,\omega^{2/\alpha}). $ Moreover, there
exist bounded operators $A_n$ and $\alpha$-times resolvent
families $S_{\alpha,n}(t)$ for $A_n$ satisfying $
||S_{\alpha,n}(t) || \leq MCe^{(2\omega)^{1/\alpha} t},$ for all
$t\geq 0,~n\in \mathbb{N},$ and
$$
S_{\alpha,n}(t)x \to S_{\alpha}(t)x \quad \mbox{as} \quad n\to
+\infty
$$
for all $x \in B,\; t\geq 0.$ Moreover, the convergence is uniform
in $t$ on every compact subset of $ \mathbb{R}_+$.
\end{theor}

In the following, we denote by $\Sigma_{\theta}(\omega)$ the open
sector with vertex $\omega \in \mathbb{R}$ and opening angle
$2\theta$ in the complex plane which is symmetric with respect to
the real positive axis, i.e.
$$ \Sigma_{\theta}(\omega) := \{ \lambda \in \mathbb{C}:
|arg(\lambda -\omega)| < \theta \}.$$

We recall from \cite[Definition 2.13]{Ba} that an $\alpha$-times
resolvent family $S_{\alpha}(t)$ is called {\tt analytic} if
$S_{\alpha}(t)$ admits an analytic extension to a sector
$\Sigma_{\theta_0}$ for some $\theta_0 \in (0, \pi/2].$ An
$\alpha$-times analytic resolvent family is said to be of {\tt
analyticity type} $(\theta_0, \omega_0)$ if for each $\theta <
\theta_0$ and $\omega > \omega_0$ there is $M= M(\theta, \omega)$
such that
$$ \| S_{\alpha}(t)\| \leq Me^{\omega Re t}, \quad t \in
\Sigma_{\theta}.$$ The set of all operators $A \in
\mathcal{C}^{\alpha}$ generating $\alpha$-times analytic resolvent
families $S_{\alpha}(t)$ of type $(\theta_0, \omega_0)$ is denoted
by $ \mathcal{A}^{\alpha}(\theta_0, \omega_0).$ In addition,
denote $ \mathcal{A}^{\alpha}(\theta_0):= \bigcup \{
\mathcal{A}^{\alpha}(\theta_0, \omega_0); \omega_0 \in
\mathbb{R}_+ \}, \quad \mathcal{A}^{\alpha} := \bigcup \{
\mathcal{A}^{\alpha}(\theta_0); \theta_0 \in (0, \pi/2] \}.$ For
$\alpha =1$ we obtain the set of all generators of analytic
semigroups.

\begin{rem}\label{rem3}{\rm We note that the spatial regularity condition
$\mathcal{R}(S_\alpha(t)) \subset D(A)$ for all $t>0$
is satisfied by $\alpha$-times resolvent families whose generator
$A$ belongs to the set $\mathcal{A}^{\alpha}(\theta_0, \omega_0)$
where $ 0<\alpha <2$ (see \cite[ Proposition 2.15]{Ba}). In
particular, setting $\omega_0=0$ we have that $ A\in
\mathcal{A}^{\alpha}(\theta_0, 0)$ if and only if $-A$ is a
positive operator with spectral angle less or equal to $ \pi -
\alpha(\pi/2+\theta).$ Note that such condition is also equivalent
to the following
\begin{equation}{\label{eq29}}
 \Sigma_{\alpha(\pi/2 +\theta)} \subset \rho(A) \mbox{ and } \|
 \lambda (\lambda I - A)^{-1} \| \leq M, \quad \lambda \in
 \Sigma_{\alpha(\pi/2 +\theta)}.
\end{equation}
}
\end{rem}

The above considerations give us the following remarkable
corollary.

\begin{cor} \label{cor5}
Suppose $A$ generates an analytic semigroup of angle $\pi/2$ and
$\alpha \in (0,1)$. Then $A$ generates an $\alpha$-times analytic
resolvent family.
\end{cor}

\bgproof Since $A$ generates an analytic semigroup of angle
$\pi/2$ we have
$$ \|
 \lambda (\lambda I - A)^{-1} \| \leq M, \quad \lambda \in
 \Sigma_{\pi -\epsilon}.$$

Then the condition (\ref{eq29}) (see also \cite[Corollary
2.16]{Ba}) implies $A \in \mathcal{A}^{\alpha}(\min
\{\frac{2-\alpha}{2\alpha}\pi, \frac{1}{2}\pi \},0)$, $\alpha \in
(0,2), $ that is $A$ generates an $\alpha$-times analytic
resolvent family.
\edproof

In the sequel we will use the following assumptions concerning
Volterra equations:
\begin{description}
\item[(A1)] $A$ is the generator of $C_0$-semigroup
 and $\alpha\in (0,1)$; \hspace{2ex} or
\item[(A2)] $A$ is the generator of a strongly continuous
 cosine family and $\alpha\in (0,2)$.
\end{description}

Observe that (A2) implies (A1) but not vice versa.

\section{Weak vs.\ mild solutions}\label{ScylW}

Assume that $H$ and $U$ are separable Hilbert spaces. Let the
cylindrical Wiener process $W$ be defined on a stochastic basis
$(\Omega,\mathcal{F},(\mathcal{F})_{t\geq 0},P)$, with the
positive symmetric covariance operator $Q\in L(U)$. This is known
that the process $W$ takes values in some superspace of $U$. (For
more details concerning cylindrical Wiener process we refer to
\cite{DPZa} or \cite{Ka}.)

We define the subspace $U_0:=Q^{1/2}(U)$ of the space $U$,
 endowed with the inner
product $\langle u,v\rangle_{U_0}:=\langle Q^{-1/2}u,
Q^{-1/2}v\rangle_U$. The set $L_2^0:=L_2(U_0,H)$ of all
Hilbert-Schmidt operators from $U_0$ into $H$, equipped with the
norm $|C|_{L_2(U_0,H)}:= (\sum_{k=1}^{+\infty} |Cf_k|_H^2)^{1/2}$,
where $\{f_k\}\subset U_0$ is an orthonormal basis of $U_0$,
is a separable Hilbert space. We assume that $\Psi$ belongs to the
class of measurable $L_2^0$-valued processes.

By $\mathcal{N}^2(0,T;L_2^0)$ we denote a Hilbert space of all
$L_2^0$-predictable processes $\Psi$ such that $|| \Psi||_T <
+\infty$, where
\begin{eqnarray*}
||\Psi||_T &:=& \left\{\mathbb{E}\left( \int_0^T
|\Psi(\tau)|_{L_2^0}^2\,d\tau \right) \right\}^{\frac{1}{2}} \\
& = & \left\{\mathbb{E} \int_0^T \left[ \mathrm{Tr}
(\Psi(\tau)Q^{\frac{1}{2}}) (\Psi(\tau)Q^{\frac{1}{2}})^* \right]
d\tau \right\}^{\frac{1}{2}}.
\end{eqnarray*}

\noindent We shall use the following
 {\sc Probability Assumptions} (abbr. (PA)):
\begin{enumerate}
\item $X_0$ is an $H$-valued, $\mathcal{F}_0$-measurable random variable;
\item $\Psi\in \mathcal{N}^2(0,T;L_2^0)$ and the interval $[0,T]$ is fixed.
\end{enumerate}

\begin{defin} \label{dSW4}
Assume that (PA) hold. An $H$-valued predictable process
$X(t),~t\in [0,T]$, is said to be a ~{\tt strong solution}~ to
(\ref{eq4}), if $X$ takes values in $D(A)$, $P$-a.s.,
\begin{equation} \label{eSW3.1}
\mbox{for~any~} t\in [0,T], \quad
\int_0^t |g_\alpha(t-\tau)AX(\tau)|_H \,d\tau<+\infty,\quad
P-a.s., \quad \alpha >0,
\end{equation}
and for any $t\in [0,T]$ the equation (\ref{eq4}) holds $P$-a.s.
\end{defin}

Let $A^*$ denote the adjoint of $A$ with a dense domain
$D(A^*)\subset H$ and the graph norm $|\cdot |_{D(A^*)}$.

\begin{defin} \label{dSW5}
Let (PA) hold. An $H$-valued predictable process $X(t),~t\in
[0,T]$, is said to be a {\tt weak solution} to (\ref{eq4}), if
$P(\int_0^t|g_\alpha(t-\tau)X(\tau)|_H d\tau<+\infty)=1,~ \alpha
>0$, and if for all $\xi\in D(A^*)$ and all $t\in [0,T]$ the
following equation holds
$$
\langle X(t),\xi\rangle_H = \langle X_0,\xi\rangle_H + \langle
\int_0^t g_\alpha(t-\tau)X(\tau)\,d\tau, A^*\xi\rangle_H + \langle
\int_0^t \Psi(\tau)dW(\tau),\xi\rangle_H, ~~P\mathrm{-a.s.}
$$
\end{defin}

\begin{defin} \label{dSW6}
Assume that $X_0$ is $\mathcal{F}_0$-measurable random variable.
An $H$-valued predictable process
$X(t),~t\in [0,T]$, is said to be a {\tt mild solution} to the
stochastic Volterra equation (\ref{eq4}), if~ $ \mathbb{E}(
\int_0^t |S_\alpha(t-\tau) \Psi(\tau)|_{L_2^0}^2 \,d\tau
)<+\infty, ~ \alpha >0$, for $t\leq T$ and, for arbitrary $t\in
[0,T]$,
\begin{equation}\label{eSW9}
X(t) = S_\alpha(t)X_0 + \int_0^t
S_\alpha(t-\tau)\Psi(\tau)\,dW(\tau), \quad P-a.s.
\end{equation}
where $S_\alpha(t)$ is the $\alpha$-times resolvent family.
\end{defin}

We will use the following result.

\begin{prop} \label{pSW1} (see, e.g.\cite[Proposition 4.15]{DPZa})\\
Assume that $A$ is a closed linear unbounded operator with the dense domain
$D(A)\subset H$. Let $\Phi(t), t\in [0,T]$, be an $L_2(U_0,H)$-predictable
process.
If $~\Phi(t)\in D(A), ~~P-a.s.$ for all $t\in [0,T]$ and
$$ P\left( \int_0^T | \Phi(s)|_{L_2^0}^2\,ds <\infty \right) =1,~~
P\left( \int_0^T |A \Phi(s)|_{L_2^0}^2\,ds <\infty \right) =1,
$$
then $\quad \displaystyle P\left( \int_0^T \Phi(s)\,dW(s) \in D(A)
\right) =1~~$ and
$$ A \int_0^T \Phi(s)\,dW(s) = \int_0^T A\Phi(s)\,dW(s), \quad P-a.s.
$$
\end{prop}

We define the stochastic convolution
\begin{equation} \label{eSW18a}
W_\alpha^\Psi(t) := \int_0^t S_\alpha(t-\tau)\Psi(\tau)\,dW(\tau),
\end{equation}
where $\Psi\in\mathcal{N}^2(0,T;L_2^0)$. Because $\alpha$-times
resolvent families $S_\alpha(t),~t\geq 0$, are bounded, then
$S_\alpha(t-\cdot)\Psi (\cdot)\in \mathcal{N}^2(0,T;L_2^0)$, too.

Analogously like in \cite{Ka}, we can formulate the following result.

\begin{prop} \label{pro2}
 Assume that $S_\alpha (t), t\ge 0$, are the resolvent operators to
 (\ref{eq3a}). Then, for any process $\Psi \in \mathcal{N}^2(0,T;L_2^0)$,
 the convolution $W_\alpha^\Psi (t), ~t\ge 0, ~\alpha >0$, given by
 (\ref{eSW18a}) has a predictable version. Additionally, the process
 $W_\alpha^\Psi (t), ~t\ge 0, ~\alpha >0$, has square integrable trajectories.
\end{prop}

Under some conditions a mild solution to Volterra equations is a
weak solution and vice versa, see \cite[Propositions 4 and 5]{Ka}.

Now, we can prove that a mild solution to the equation~(\ref{eq4})
is a weak solution to~(\ref{eq4}).
\begin{prop} \label{pSW4}
If $\Psi\in\mathcal{N}^2(0,T;L_2^0)$ and
$\Psi(\cdot,\cdot)(U_0)\subset D(A)$, $P$-a.s., then the
stochastic convolution $W_\alpha^\Psi (t),~t\ge 0, ~\alpha >0$, given
by  (\ref{eSW18a}), fulfills the equation
\begin{equation}\label{eSW18b}
 \langle W_\alpha^\Psi(t),\xi\rangle_H =
\int_0^t \langle g_\alpha(t-\tau)W_\alpha^\Psi(\tau),
A^*\xi\rangle_H + \int_0^t \langle
\xi,\Psi(\tau)dW(\tau)\rangle_H, \quad \alpha\in (0,2),
\end{equation}
for any $t\in [0,T]$ and $\xi\in D(A^*)$.
\end{prop}
\bgproof Let us notice that the process $W_\alpha^\Psi$ has integrable trajectories.
For any $\xi\in D(A^*)$ we have
\begin{eqnarray*}
 \int_0^t \langle g_\alpha(t-\tau)W_\alpha^\Psi(\tau),A^*\xi\rangle_Hd\tau \!&\!\equiv\!&\!
 ~\mbox{(from~ (\ref{eSW18a}))} \\
 \!&\!\equiv\!&\! \int_0^t \langle g_\alpha(t-\tau) \int_0^\tau S_\alpha(\tau-\sigma)\Psi(\sigma)
 dW(\sigma),A^*\xi\rangle_Hd\tau =\\
  (\mbox{from~Dirichlet's~formula}\!&\! \mbox{and} \!&\!
  \mbox{stochastic Fubini's theorem})  \\
 \!&\!=\!&\! \int_0^t \langle\left[\int_\sigma^t g_\alpha(t-\tau)S_\alpha(\tau-\sigma)d\tau\right]
 \Psi(\sigma)  dW(\sigma),A^*\xi\rangle_H\\
 \!&\!=\!&\! \langle \!\int_0^t\! \left[\int_0^{t-\sigma}\! g_\alpha(t-\sigma-z)S_\alpha(z)dz\right]
 \!\Psi(\sigma) dW(\sigma),A^*\xi\rangle_H\\
 (\mbox{where~} z:=\tau-\sigma\!&\!\!&\! \mbox{~and~from definition~of~convolution})\\
  \!&\!=\!&\! \langle \int_0^t A [ (g_\alpha\star S_\alpha)(t-\sigma)]\Psi(\sigma) dW(\sigma),
  \xi\rangle_H = \\
  (\mbox{from~the~resolvent~equation~(\ref{eq4a})} \!&\!\!&\!  \mbox{because~~}
     A(g_\alpha\star S_\alpha)(t-\sigma)x = (S_\alpha(t-\sigma)-I)x, \\
    \!&\!\!&\! ~\mbox{~where~} x\in D(A))  \\
 \!&\!=\!&\! \langle \int_0^t [S_\alpha(t-\sigma)-I]\Psi(\sigma) dW(\sigma),\xi\rangle_H = \\
 = \langle \int_0^t S_\alpha(t-\sigma)\Psi(\sigma) dW(\sigma),\xi\rangle_H
 \!&\!-\!&\!\langle \int_0^t \Psi(\sigma) dW(\sigma),\xi\rangle_H .
\end{eqnarray*}
Hence, we obtained the following equation
$$
 \langle W_\alpha^\Psi(t),\xi\rangle_H = \int_0^t \langle g_\alpha(t-\tau)W_\alpha^\Psi(\tau),
 A^*\xi\rangle_Hd\tau + \int_0^t \langle \xi,\Psi(\tau)dW(\tau)\rangle_H
$$
for any $\xi\in D(A^*)$. \edproof

Immediately from the equation (\ref{eSW18b}) we deduce the following result.
\begin{cor} \label{cor2a}
 If $A$ is a bounded operator and $\Psi\in\mathcal{N}^2(0,T;L_2^0)$, then the
 following equality holds
 \begin{equation} \label{eq21a}
 W_\alpha^\Psi(t) =
 \int_0^t g_\alpha(t-\tau) A W_\alpha^\Psi(\tau)d\tau +
 \int_0^t\Psi(\tau)dW(\tau),
\end {equation}
for $t\in [0,T], ~\alpha >0$.
\end{cor}

\begin{rem}{\em The formula (\ref{eq21a}) says that the convolution
$W_\alpha^\Psi(t), ~t\ge 0, ~\alpha >0$, is a strong solution to
(\ref{eq4}) with $X_0\equiv 0$ if the operator $A$ is bounded.}
\end{rem}

\section{Strong solutions}

In this section we provide sufficient conditions under which the
stochastic convolution $W_\alpha^\Psi (t)$, $t\ge 0$, $\alpha >0$,
defined by (\ref{eSW18a}) is a strong solution to the equation
(\ref{eq4}).

\begin{lem} \label{pSW5}
Let $A$ be a closed linear unbounded operator with dense domain $D(A)$
equipped with the graph norm $|\cdot|_{D(A)}$.
Assume that {\bf (A1)} or {\bf (A2)} holds. If $\Psi$ and $A\Psi$
belong to $\mathcal{N}^2(0,T;L_2^0)$ and in addition
$\Psi(\cdot,\cdot)(U_0)\subset D(A),$ P-a.s., then (\ref{eq21a}) holds.
\end{lem}
\bgproof Because formula (\ref{eq21a})
holds for any bounded operator, then it holds for the Yosida
approximation $A_n$ of the operator $A$, too, that is
$$ W_{\alpha,n}^\Psi(t) =
\int_0^t g_\alpha(t-\tau) A_n W_{\alpha,n}^\Psi(\tau)d\tau +
\int_0^t\Psi(\tau)dW(\tau), $$ where
$$ W_{\alpha,n}^\Psi(t) := \int_0^t S_{\alpha,n}(t-\tau)\Psi(\tau)dW(\tau).$$
By Proposition~\ref{pSW1}, we have
$$ A_n W_{\alpha,n}^\Psi(t) =
A_n \int_0^t S_{\alpha,n}(t-\tau)\Psi(\tau)dW(\tau).
$$
By assumption $\Psi\in \mathcal{N}^2(0,T;L_2^0)$. Because the
operators $S_{\alpha,n}(t)$ are deterministic and bounded for any
$t\in [0,T], ~\alpha >0, ~n\in\mathbb{N}$, then the operators
$S_{\alpha,n}(t-\cdot )\Psi(\cdot)$ belong to
$\mathcal{N}^2(0,T;L_2^0)$, too. In consequence, the difference
\begin{equation}\label{eq21b}
 \Phi_{\alpha,n}(t-\cdot ) := S_{\alpha,n}(t-\cdot )\Psi(\cdot)
  - S_\alpha(t-\cdot )\Psi(\cdot)
\end{equation}
belongs to $\mathcal{N}^2(0,T;L_2^0)$ for any $t\in [0,T], ~\alpha >0$ and
$n\in\mathbb{N}$. This means that
\begin{equation}\label{eq22}
 \mathbb{E}\left(\int_0^t |\Phi_{\alpha,n}(t-\tau)|_{L_2^0}^2d\tau \right)
 < +\infty
\end{equation}
for any $t\in [0,T]$.

Let us recall that the cylindrical Wiener process $W(t)$, $t\ge 0$,
can be written in the
form
\begin{equation}\label{eq23}
 W(t) =\sum_{j=1}^{+\infty} f_j\,\beta_j(t),
\end{equation}
where $\{f_j\}$ is an orthonormal basis of $U_0$ and $\beta_j(t)$ are
independent
real Wiener processes. From (\ref{eq23}) we have
\begin{equation}\label{eq24}
 \int_0^t \Phi_{\alpha,n}(t-\tau)\,dW(\tau) = \sum_{j=1}^{+\infty}
 \int_0^t \Phi_{\alpha,n}(t-\tau)\,f_j\,d\beta_j(\tau).
\end{equation}
In consequence, from (\ref{eq22})
\begin{equation}\label{eq25}
 \mathbb{E}\left[\int_0^t \left( \sum_{j=1}^{+\infty}
 |\Phi_{\alpha,n}(t-\tau)\,f_j|_H^2 \right) d\tau \right]
 < +\infty
\end{equation}
for any $t\in [0,T]$. Next, from (\ref{eq24}), properties of stochastic
integral and (\ref{eq25}) we obtain for any $t\in[0,T]$,
\begin{eqnarray*}
 \mathbb{E}\left| \int_0^t \Phi_{\alpha,n}(t-\tau)\,dW(\tau) \right|_H^2
 &=& \mathbb{E}\left| \sum_{j=1}^{+\infty}\int_0^t
 \Phi_{\alpha,n}(t-\tau)\,f_j\,d\beta_j(\tau) \right|_H^2  \le \\
  \mathbb{E}\left[ \sum_{j=1}^{+\infty} \int_0^t
 |\Phi_{\alpha,n}(t-\tau)\,f_j|_H^2 d\tau \right]
 &\le & \mathbb{E}\left[ \sum_{j=1}^{+\infty} \int_0^T
  |\Phi_{\alpha,n}(T-\tau)\,f_j|_H^2 d\tau \right] <+\infty.
\end{eqnarray*}

By Theorem \ref{th2}, the convergence (\ref{eq10}) of
$\alpha$-times resolvent families is uniform in $t$ on every
compact subset of $\mathbb{R}_+$, particularly on the interval
$[0,T]$. Now, we use (\ref{eq10}) in the Hilbert space $H$, so
(\ref{eq10}) holds for every $x\in H$. Then, for any fixed $\alpha$ and $j$,
\begin{equation}\label{eq26}
 \int_0^T |[S_{\alpha,n}(T-\tau)-S_\alpha(T-\tau)]\,
 \Psi(\tau)\,f_j|_H^2 d\tau
\end{equation}
tends to zero for $n\to +\infty$.
So, summing up our considerations, particularly using (\ref{eq25}) and
(\ref{eq26}) we can write
\begin{eqnarray*}
\sup_{t\in [0,T]} \, \mathbb{E}\left| \int_0^t
\Phi_{\alpha,n}(t-\tau)dW(\tau)
 \right|_H^2 \!\equiv \!\sup_{t\in [0,T]}\, \mathbb{E}\left| \int_0^t
 [S_{\alpha,n}(t-\tau)-S_\alpha(t-\tau)] \Psi(\tau)dW(\tau)\right|_H^2
 \!\!& \!\!\le\!\! &\!\! \\ \le
 \mathbb{E}\left[ \sum_{j=1}^{+\infty} \int_0^T
 | [ S_{\alpha,n}(T-\tau)-S_\alpha(T-\tau)]\Psi(\tau)\,f_j |_H^2
 d\tau \right] \!\!&\!\!\to \!\!& \!\!0
\end{eqnarray*}
as $n \to +\infty$ for any fixed $\alpha >0$.

Hence, by the Lebesgue dominated convergence theorem we obtained
\begin{equation}\label{eq27}
 \lim_{n\to +\infty} \sup_{t\in [0,T]} \mathbb{E} \left|
 W_{\alpha,n}^\Psi(t)- W_\alpha^\Psi(t)\right|_H^2 =0.
\end{equation}
By Proposition \ref{pSW1}, $P(W_\alpha^\Psi (t)\in D(A))=1$.

For any $n\in\mathbb{N}$, $\alpha >0$, $t\ge 0$, we have
$$ |A_n W_{\alpha,n}^\Psi (t) - A W_\alpha^\Psi (t)|_H \le
    N_{n,1}(t) +  N_{n,2}(t), $$
where
\begin{eqnarray*}
 N_{n,1}(t) & := & |A_n W_{\alpha,n}^\Psi (t) - A_n W_\alpha^\Psi (t)|_H , \\
 N_{n,2}(t) & := & |A_n W_\alpha^\Psi (t) - A W_\alpha^\Psi (t)|_H =
             |(A_n-A)W_\alpha^\Psi (t)|_H  \,.
\end{eqnarray*}
Then
\begin{equation}\label{eq28}
 |A_n W_{\alpha,n}^\Psi (t) - A W_\alpha^\Psi (t)|_H^2 \le
 N_{n,1}^2 (t) + 2 N_{n,1}(t) N_{n,2}(t) + N_{n,2}^2(t) .
\end{equation}

Let us study the term $N_{n,1}(t)$. Note that, either in cases
{\bf (A1)} or {\bf (A2)} the unbounded operator $A$ generates a
semigroup. Then we have for the Yosida approximation the following
properties:
\begin{equation}\label{eq30}
 A_nx=J_nAx \quad \mbox{for~any~} x\in D(A), \quad \sup_n ||J_n||
 < \infty
\end{equation}
where $A_nx=nAR(n,A)x=AJ_nx$ for any $x\in H$, with
$J_n:=nR(n,A).$ Moreover (see \cite[Chapter II, Lemma
3.4]{En-Na}):
\begin{eqnarray}
 \lim_{n\to\infty}J_nx &=& x \qquad \mbox{for~any~} x\in H, \nonumber \\
 \lim_{n\to\infty} A_nx &=& Ax \qquad \mbox{for~any~} x\in D(A).
 \label{eq31}
\end{eqnarray}
Note that $AS_{\alpha,n}(t)x = S_{\alpha,n}(t)Ax$ for all $x \in
D(A),$  since $e^{t A_n}$ commutes with $A$ and $A$ is closed (see
(\ref{eq14})). So, by Proposition \ref{pSW1} and again the
closedness of $A$  we can write
\begin{eqnarray*}
 A_n W_{\alpha,n}^\Psi(t) &\equiv& A_n
 \int_0^t S_{\alpha,n}(t-\tau)\Psi(\tau)dW(\tau) \\
 &=& J_n \int_0^t AS_{\alpha,n}(t-\tau)\Psi(\tau)dW(\tau)
 = J_n \left[\int_0^t S_{\alpha,n}(t-\tau)A\Psi(\tau)dW(\tau)\right].
\end{eqnarray*}
Analogously,
$$  A_n W_\alpha^\Psi(t) =
  J_n \left[\int_0^t S_\alpha(t-\tau)A\Psi(\tau)dW(\tau)\right].
$$
By (\ref{eq30}) we have
\begin{eqnarray*}
N_{n,1}(t) &=& |J_n\int_0^t [S_{\alpha,n}(t-\tau)-S_\alpha(t-\tau)]
  A\Psi(\tau)dW(\tau)|_H \\
  &\le & |\int_0^t [S_{\alpha,n}(t-\tau)-S_\alpha(t-\tau)]
  A\Psi(\tau)dW(\tau)|_H \;.
\end{eqnarray*}
From assumptions, $A\Psi \in \mathcal{N}^2(0,T;L_2^0)$. Then the term
$[S_{\alpha,n}(t-\tau)-S_\alpha(t-\tau)]A\Psi(\tau)$ may be estimated
like the difference $\Phi_{\alpha,n}$ defined by (\ref{eq21b}).

Hence, from (\ref{eq30}) and (\ref{eq27}),
for the first term of the right hand side of (\ref{eq28}) we obtain
$$ \lim_{n\to +\infty}\;\; \sup_{t\in [0,T]}
 \mathbb{E}(N_{n,1}^2 (t)) \to 0. $$
For the second and third terms of (\ref{eq28})
we can follow the same steps as above for proving (\ref{eq27}).
We have to use the properties of Yosida approximation, particularly the
convergence (\ref{eq31}).
So, we can deduce that
$$ \lim_{n\to +\infty} \;\;\sup_{t\in [0,T]} \mathbb{E}
  |A_n W_{\alpha,n}^\Psi(t)-AW_\alpha^\Psi(t)|_H^2=0 , $$
what gives (\ref{eSW18a}). \hfill $\blacksquare$

Now, we are able to formulate the main result of this section.

\begin{theor} \label{coSW4}
Suppose that assumptions of Lemma \ref{pSW5} hold. Then the
equation (\ref{eq4}) with $X_0\equiv 0$ has a strong solution.
Precisely, the convolution $W_\alpha^\Psi$ defined by
(\ref{eSW18a}) is the strong solution to~(\ref{eq4}) with $X_0\equiv 0$.
\end{theor}
\bgproof We have to show only the condition (\ref{eSW3.1}).
By Proposition~\ref{pro2}, the convolution
$W_\alpha^\Psi (t)$, $t\ge 0$, $\alpha >0$, has integrable trajectories,
that is, $W_\alpha^\Psi (\cdot )\in L^1([0,T];H)$,
P-a.s. The closed linear unbounded operator $A$ becomes bounded on
($D(A),|\cdot|_{D(A)}$), see \cite[Chapter 5]{We}. So, we obtain
$AW_\alpha^\Psi (\cdot )\in L^1([0,T];H)$, P-a.s. Hence, the
function $g_\alpha(T-\tau)AW_\alpha^\Psi (\tau)$ is integrable
with respect to $\tau$, what finishes the proof. \edproof

The following result is an immediate consequence of Corollary
\ref{cor5}
 and Theorem \ref{coSW4}.

\begin{cor} \label{cor6}
Assume that $A$ generates an analytic semigroup of angle $\pi/2$
and $\alpha \in (0,1)$.
 If $X_0=0$, then the equation (\ref{eq4}) has a strong solution.
\end{cor}

{\bf Acknowledgement} The authors would like to thank the
referee for the careful reading of the manuscript.
The valuable remarks made numerous improvements throughout.


\begin{thebibliography}{99}

\bibitem{Ba} E.~Bazhlekova, {\it Fractional Evolution Equations in
Banach Spaces,} Ph.D.~Dissertation, Eindhoven University of
Technology, 2001.


\bibitem{Bo-Tu} S.~Bonaccorsi, L.~Tubaro, {\it Mittag-Leffler's
function and stochastic Volterra equations of convolution type,}
Stochastic Anal.\ Appl.\ \textbf{21} (1) (2003), 61-78.

\bibitem{Cl-Pr96} Ph.~Cl\'{e}ment, G.~Da Prato, {\it Some results on
  stochastic convolutions arising in Volterra equations perturbed by noise},
Rend.\ Math.\ Acc.\ Lincei s.\ 9, \textbf{7} (1996), 147-153.

\bibitem{DPZa} G.~Da Prato, J.~Zabczyk, \textit{Stochastic
equations in infinite dimensions}, Cambridge University Press,
Cambridge, 1992.

\bibitem{En-Na} K.J.~Engel, R.Nagel, \textit{One-parameter semigroups for linear evolution equations}, Graduate texts in Mathematics
{\bf 194}, Springer, New York, 2000.

\bibitem{Fu} Y.~Fujita, {\it Integrodifferential equations which
 interpolates the heat equation and the wave equation},
 J.~Math.\ Phys.\ \textbf{30} (1989), 134--144.

 \bibitem{Ho-Ok-Ub-Zh} H.~Holden, B.~$\O$ksendal, J. Ub$\o$e, T. Zhang,
{\it Stochastic Partial Differential Equations: A modeling, white
noise functional approach,} Probability and its applicatons,
Birkh\"auser, 1996.

\bibitem{Ka} A.~Karczewska, {\it Properties of convolutions arising in
stochastic Volterra equations}, preprint:
http://xxx.lanl.gov/ps/math.PR/0509012

\bibitem{KL2} A.~Karczewska, C.~Lizama, {\it Strong
solutions to stochastic Volterra equations}, submitted.

\bibitem{LiZh} M.~Li, Q.~Zheng, {\it On spectral inclusions and
approximations of $\alpha$-times resolvent families,} Semigroup
Forum {\bf 69} (2004), 356-368.

\bibitem{Pa} A. Pazy, \textit{Semigroups of Linear Operators
and Applications to Partial Differential Equations}, Springer, New
York, 1983.

\bibitem{We} J.~Weidmann, \textit{Linear operators in Hilbert spaces},
Springer, New York, 1980.


\end{thebibliography}
\end{document}